\documentclass[12pt]{amsart}
\usepackage{amsmath,amsthm,epsfig,amssymb}
\usepackage{amsthm,amsfonts,amssymb,amscd}
\usepackage[all]{xy}
\usepackage{young}
\newtheorem{thm}[subsection]{Theorem}
\newtheorem{prop}[subsection]{Proposition}

\newtheorem{lem}[subsection]{Lemma}
\newtheorem{corol}[subsection]{Corollary}
\newtheorem{rem}[subsection]{Remark}
\theoremstyle{definition}
\newtheorem{Def}[subsection]{Definition}
\newtheorem{Not}[subsection]{Notation}

\newtheorem{proposition-definition}[subsection]{Proposition-Definition}

\begin{normalsize} \end{normalsize}

\newcommand{\CC}{{\mathbb C}}

\newcommand{\PP}{{\mathbb P}}
\newcommand{\NN}{{\mathbb N}}

\newcommand{\OOO}{{\mathcal O}}

\newcommand{\AAA}{{\mathcal L}}
\newcommand{\WWW}{{\mathcal L}}

\newcommand{\GGG}{{\mathcal G}}

\numberwithin{equation}{section}
\renewcommand\square{\frame{\phantom{{\large x}}}}

\author{F. Laytimi}

\address{F. L.: Math\'ematiques - b\^{a}t. M2, Universit\'e Lille 1,
F-59655 Villeneuve d'Ascq Cedex, France}
\email{fatima.laytimi@math.univ-lille1.fr}

\author{W. Nahm}

\address{W. N.: Dublin Institute for Advanced Studies,
10 Burlington Road, Dublin 4, Ireland} \email{wnahm@stp.dias.ie}


\title{ Semiample and $k$-ample vector bundles}

\begin{document}

\date{}

\begin{abstract}
We show that tensor products of semiample vector bundles are semiample. 
For $k$-ampleness in the sense of Sommese, we show that over compact complex manifolds 
tensor products of semiample and $k$-ample vector bundles are $k$-ample, 
and the sum of $k$-ample vector bundles is $k$-ample. In particular,  results of 
Sommese on $k$-ampleness are strenghtened. 
\end{abstract}

\maketitle

\section{Tensor products of semiample vector bundles} \setcounter{page}{1}

\begin{Def}
Let $X$ be a complex space. A line bundle $L$ on $X$ is said to be 
semiample if for some $r>0, \, L^r $ is generated by sections.
A vector bundle $E$ is  semiample if $\OOO_{\PP E}(1)$ is semiample. 
\end{Def}

\begin{Not}\label{Notgbs}
Since we often need to use the phrase ``generated by sections'', we will abbreviate it by
``gbs''.
\end{Not}

We start  by giving several interpretations of the statement:
$$ \OOO_{\PP E}(r)  \textrm {\ is gbs}.$$ 

1) \  Let $\pi:\PP E\rightarrow X$ be the projection. The fibre of $\pi$ over a point 
$x \in X$ is given by the 1-codimensional
subspaces $V$ of $E_x$, where $E_x$ is the fibre of $E$ over $x$. Points of $\PP E$ will
be denoted by pairs $(x,V)$.
Let $V=v^\perp$ for $v \in E_x^*$, and let ${\tilde v} $ be the corresponding element of $(E_x/V)^*$.  A section $ s\in H^0(\PP E,\OOO_{\PP E}(r))$ generates 
the fibre of $\OOO_{\PP E}(r)$ 
over $(x,V)$ iff
$$\langle s(x,V), {\tilde v}^r\rangle\neq 0,$$
where $\langle\  , \ \rangle$ denotes the pairing of dual vector spaces.\\

2)\  Equivalently, let $\tilde s$ be the section in $H^0(X,S^rE)$ which corresponds to
$s$ by the natural isomorphism
$$H^0(\PP E,\OOO_{\PP E}(r))\simeq H^0(X,S^rE).$$ 
Then $s$ generates the fibre of $\OOO_{\PP E}(r)$ over $(x,v^\perp)$ iff
$$\langle \tilde s(x), v^{\otimes r}\rangle\neq 0.$$

3)\  Given $x\in X$, the value $\tilde s(x)$ of each section 
$\tilde s \in H^0(X,S^rE)$ \\
can be regarded as a homogeneous 
polynomial map of degree $r$ from $E_x^*$ to $\CC$.\\ 
If the fibre of $\OOO_{\PP E}(r)$ 
over $v^\perp$ is not generated by 
any global section, then these polynomials have a common zero at 
$v$. If $\OOO_{\PP E}(r)$
is gbs, 
then 0 is the only common zero.\\

4)\ Since $v^{\otimes r}$ annihilates the kernel of the natural map
$E^{\otimes r}\twoheadrightarrow S^rE$,
the fibre of $\OOO_{\PP E}(r)$ over $(x,v^\perp)$ is generated by a section iff
there is a section $\tilde s \in H^0(X,E^{\otimes r})$ such that
$$\langle \tilde s(x), v^{\otimes r}\rangle\neq 0.$$

If $\OOO_{\PP E}(r)$ is gbs, then $\OOO_{\PP E}(kr)$ is gbs for every positive 
integer $k$. All these facts will be used in the following without reference.

\begin{thm} Direct sum and tensor product of semiample vectors bundles
on a complex space $X$ are semiample.
\end{thm}
This result is implied by the following two lemmata.

\begin{lem} \label{sumgbs} Let $X$ be a complex space and $E,F$ vector bundles on $X$. 
Suppose that $\OOO_{\PP E}(r)$ and $\OOO_{\PP F}(r)$ are gbs \rm (see Notation \ref{Notgbs}\rm )
for some $r>0$. Then $\OOO_{\PP(E\oplus F)}(r)$ is gbs. 
\end{lem}

{\it Proof}:
For each $x\in X$ and each non-zero $u\in (E_x^*\oplus F_x^*)$ we have to find a
section
in $H^0(X,S^r(E\oplus F))$ which is not annihilated by $u^{\otimes r}$. Consider 
$u=u_E+u_F$,
where $u_E,u_F$ are the components in $E_x^*, F_x^*$ resp. We may assume that 
$u_E\neq 0$. Since $\OOO_{\PP E}(r)$ 
is gbs, there is a section in $H^0(X,S^rE)$ which is not annihilated by $u_E^{\otimes r}$. 
The image of this section 
under the natural injection $H^0(X, S^rE)\hookrightarrow H^0(X, S^r(E\oplus F))$ is not 
annihilated by $u^{\otimes r}$.
$\hfill{\square}$

\begin{lem} \label{gbs} Let $X$ be a complex space and $E,F$ vector bundles on $X$. 
Suppose that  $\OOO_{\PP E}(r)$ and $\OOO_{\PP F}(r)$ are gbs \rm (see Notation \ref{Notgbs}\rm ) for some $r>0$. 
Then there exists a positive integer $n$ depending only on $r$ and $e={\rm rk} E$ such that 
 $\OOO_{\PP (E\otimes F)}(nr)$ is gbs.
\end{lem}

The proof of this lemma needs some preparations. We will use the terminology
of decorated oriented graphs (see Feynman diagrams in physics). Decoration means 
that to each vertex and to each arrow another object is associated. The set of 
arrows of a graph $\gamma$ will be denoted by $A(\gamma)$. 
Any arrow $a$ is said to have a head $h(a)$ and a tail $t(a)$. A vertex which is the
head of exactly $n$ arrows is said to have indegree $n$, a vertex which is the tail 
of exactly $n$ arrows is said to have outdegree $n$.
The bidegree of a vertex is written as (indegree, outdegree).

\begin{Def} 
${}$

Let $V,W$ be vector spaces. Let
$i=1,\ldots,n$, $\tilde v_i\in S^r V$, $\tilde w_i\in S^r W$, $u\in (V\otimes W)^*$. Then
$\Gamma_{V,W}(n,r,\tilde v_1,\ldots,\tilde v_n,\tilde w_1,\ldots,\tilde w_n,u)$ is the set 
of all decorated oriented graphs 
with $2n$ fixed vertices $\alpha_i, \beta_i,$ such that for each $i$ the vertex $\alpha_i$ 
has bidegree $(0,r)$
and is decorated by $\tilde v(\alpha_i)=\tilde v_i$, the vertex $\beta_i$ has bidegree $(r,0)$
and is decorated by 
$\tilde w(\beta_i)=\tilde w_i$, 
and such that all $rn$ arrows are decorated by $u$. We will omit the indices $V,W$ of 
$\Gamma_{V,W}$, when no
ambiguity arises.
\end{Def}

\begin{Def} 
 Let $\gamma\in \Gamma(n,1,v_1,\ldots,v_n,w_1,\ldots,w_n,u)$. The value $\vert\gamma\vert$ 
 of $\gamma$ is defined by 
$$\vert\gamma\vert = \prod_{a\in A(\gamma)} u\big(v(h(a))\otimes w(t(a))\big).$$ 
\end{Def}

\begin{Def} 
 Let $\gamma\in \Gamma(n,r,\tilde v_1,\ldots,\tilde v_n,\tilde w_1,\ldots,\tilde w_n,u)$. 
 Then the value $\vert\gamma\vert$ of 
$\gamma$ is defined by the following two conditions:
Firstly, $\vert\gamma\vert$ depends linearly on all $\tilde v_i$ and $\tilde w_i$. Secondly, 
if
$\tilde v_i = v_i^{\otimes r}$, $\tilde w_i = w_i^{\otimes r}$ for $i=1,\ldots,n$, $v_i\in V$,  $w_i\in W$ 
let the expanded 
graph
$$\gamma_{ex}\in \Gamma(rn,1,\underbrace{v_1,\ldots,v_1}_{r\ times},\ldots,\underbrace{v_n,\ldots,v_n}_{r\ times},
\underbrace{w_1,\ldots,w_1}_{r\ times},\ldots,\underbrace{w_n,\ldots,w_n}_{r\ times},u)$$
be a decorated oriented graph with vertices $\alpha_i^l, \beta_i^l,$ $i=1,\ldots,n$, $l=1,\ldots,r,$ 
such that
there is a bijection $\xi:A(\gamma_{ex})\rightarrow A(\gamma)$ with
$h(a)=\alpha_i$ if $h(\xi(a))=\alpha_i^l$, $t(a)=\beta_i$ if $t(\xi(a))=\beta_i^l$.
Moreover, let
$\alpha_i^l$ be decorated by $v_i$, $\beta_i^l$ by $w_i$, for $i=1,\ldots,n$, $l=1,\ldots,r$. Then 
$$\vert\gamma\vert = \vert\gamma_{ex}\vert.$$ 
\end{Def}

For later use, note that $\gamma_{ex}$ has a symmetry group of order\\
$(r!)^{2n}s_\gamma^{-1}$,
where $s_\gamma$ is 
the order of the group of vertex preserving symmetries of $\gamma$.

\begin{prop} \label{nu} 
There is a function $\nu:\NN \times \NN \rightarrow  \NN$ 
with the following property.
Consider finite dimensional vector spaces $V,W$ with $d={\rm dim} V$ and subspaces $A\subset S^rV$,
$B \subset S^rW$, such that the corresponding spaces of polynomial maps have 0 as only 
common zero. 
Let $n=\nu(r,d)$. Then for any non-zero $u \in (V\otimes W)^*$ there is a positive
integer $m$ with $m\leq n$,
elements $\tilde v_i\in A$, $\tilde w_i\in B$ for $i=1,\ldots,m$ and a decorated directed graph
$\gamma\in \Gamma_{V,W}(m,r,\tilde v_1,\ldots,\tilde v_m,\tilde w_1,\ldots,\tilde w_m,u)$
such that  $\vert\gamma\vert\neq 0$.
\end{prop}

{\it Proof}:
We will construct a function  $\nu$ which is far from optimal, but  sufficient for our purpose.

We first show that it suffices to prove the proposition for the case that the map $\hat u\in Hom(V,W^*)$ 
corresponding to 
$u$ is bijective. Otherwise, let $V'=V/ker(\hat u)$,  $W'=W/ker(\hat u^*)$, let $\pi_V,\pi_W$ be the 
corresponding projections 
$\pi_V:S^rV\twoheadrightarrow S^rV'$, $\pi_W:S^rW\twoheadrightarrow S^rW'$ and let $u'$ be the element
of $(V'\otimes W')^*$
induced by $u$. If $$\gamma\in \Gamma(m,r,V_1,\ldots,V_m,W_1,\ldots,W_m,u), \rm {and} $$ 
$$\gamma'\in \Gamma(m,r,\pi_V(V_1),\ldots,\pi_V(V_m),\pi_W(W_1),\ldots,\pi_W(W_m),u')$$ 
have the same underlying undecorated graph, then $\vert\gamma'\vert=\vert\gamma\vert$.

If $u$ is bijective, we can identify $W$ with $V^*$ and write
$$u(v\otimes w)=\langle v,w\rangle$$  
for $v\in V$, $w\in W$. Let $SV$ be the symmetric algebra over $V$.
There are natural multiplication maps  $$m^r:S^rV\rightarrow End(SV),$$ 
and contraction maps  $$i^r:S^rW \rightarrow End(SV)$$
which restrict to
$$m^r_N:S^rV\rightarrow Hom(S^{N-r}V,S^NV),$$ 
$$i^r_N:S^rW \rightarrow Hom(S^NV,S^{N-r}V),$$ 
with $S^NV=0$ for negative $N.$ For $r=1$ the maps can be characterized by multilinearity and 
the properties
$m^1(v')v^{\otimes N}=\pi_S (v'\otimes v^{\otimes N})$, $i^1(w)v^{\otimes N}=N\langle w,v\rangle v^{\otimes (N-1)}$
for $v',v\in V$, $w\in W$, where $\pi_S$ is the projection of the tensor algebra of $V$ to the 
symmetric algebra $SV$.
For general $r$ they are characterized by multilinearity and the properties $m^r(v^{\otimes r})=m^1(v)^r$, 
$i^r(w^{\otimes r})=i^1(w)^r$.

For $r=1$ let $P$ be a product with the factors $m^1(v_i)$, $i^1(w_i)$, taken in any order, 
where $v_i\in V$, $w_i\in W$,
$i=1,\ldots,m$. For any non-negative integer $N$, the product $P$ restricts to a map $P_N\in End(S^NV)$.
For the trace of $P_N$ 
one finds
$$tr (P_N) = \sum_{\gamma\in \Gamma_m} c_\gamma\vert\gamma\vert,$$
where 
$$\Gamma_m = \Gamma_{V,W}(m,1,v_1,\ldots,v_m,w_1,\ldots,w_m,\langle,\rangle)$$
and 
$$c_\gamma={d+\rho+N-1\choose d+m-1},$$
where $\rho$ is the cardinality of the set of arrows $a$ of 
$\gamma$ such that the factor $m^1(v(h(a)))$ 
lies
to the right of the factor  $i^1(w(t(a)))$ in $P$. For $\rho>0$ this follows by induction on $\rho$
from
$i^1(w)m^1(v)=m^1(v)i^1(w)+\langle w,v\rangle$ and for $\rho=0$ from cyclic invariance of the trace and
induction on $N$.
For $n=0$ and $N=0$ the statement is obvious. The specific form of $c_\gamma$ is not important in the following.
We have given it because of its unexpected simplicity.

For arbitrary $r$ the calculation can be reduced to the case $r=1$ as in the
definition of $\vert\gamma\vert$. 
Let  $P$ be a product with the factors $m^r(\tilde v_i)$, $i^r(\tilde w_i)$, 
taken in any order, 
where $\tilde v_i\in S^rV$, $\tilde w_i\in S^rW$, $i=1,\ldots,m$. For any $N$, $P$ 
restricts to a map $P_N\in End(S^NV)$. 
To calculate the trace of $P_N$ it is sufficient to consider the case
$\tilde v_i = v_i^{\otimes r}$, 
$\tilde w_i = w_i^{\otimes r}$ for $i=1,\ldots,n$, $v_i\in V$,  $w_i\in W$. One finds
$$tr( P_N) = \sum_{\gamma\in \Gamma^r_m} c_\gamma\vert\gamma\vert,$$
where 
$$\Gamma^r_m = \Gamma_{V,W}(m,r,\tilde v_1,\ldots,\tilde v_m,\tilde w_1,
\ldots,\tilde w_m,\langle,\rangle)$$
and 
$$c_\gamma=(r!)^{2m}s_\gamma^{-1}{d+r\rho+N-1\choose d+rm-1}.$$
Here $\rho$ is the cardinality of the set of arrows $a$ of $\gamma$
such that the factor $m^r(\tilde v(h(a)))$ lies
to the right of the factor  $i^r(\tilde w(t(a)))$ in $P$, and  $s_\gamma$ is the 
order of the group of vertex preserving
symmetries of $\gamma$.

Let $\AAA^N_r$ be the subalgebra of ${\rm End}(S^NV)$ generated by the elements
$\{m^r_N(\tilde v) i^r_N(\tilde w)\vert \tilde v\in A, \tilde w\in B\} $. This 
algebra is spanned by products
$P_N=m^r_N(\tilde v_1) i^r_N(\tilde w_1)\cdots m^r_N(\tilde v_m) i^r_N(\tilde w_m)$,
where $\tilde v_i\in A$,
$\tilde w_i\in B$ for $i=1,\ldots,m$.
We have seen that the trace of $P_N$ is given by a linear combination of numbers
$\vert\gamma\vert$, where $\gamma\in \Gamma(m,r,\tilde v_1,\ldots,\tilde v_m,\tilde w_1,
\ldots,\tilde w_m,\langle,\rangle)$.
If those traces all vanish, then $ \AAA^N_r$  is a nilalgebra. In particular, there is
a non-trivial subspace $V_0\subset S^NV$
such that 
$$m^r_N(\tilde v) i^r_N(\tilde w)V_0=0$$ 
for all $\tilde v\in A$, $\tilde w\in B$.
Since the kernel of $m^r_N(\tilde v)$ vanishes for $\tilde v\neq 0$, this means that 
$i^r_N(\tilde w)V_0=0$ for all $\tilde w\in B$. 
Since the homogeneous polynomial maps in $B$ have 0 as only common zero, one can apply
a theorem of Macaulay (\cite{M} Sec. 6, or \cite{S2} p. 85 theorem 4.48). According to
this theorem the ideal 
in $SW$ generated by $B$ contains $S^NW$ if $N\geq rd$, which implies $i^N_N(\tilde w)V_0=0$
 for all $\tilde w\in S^NW$ and
yields a contradiction. Thus $\AAA^N_r$ is not nil for $N\geq rd$.

Let $\WWW^N_m$ be the subspace of $\AAA^N_r$ generated by products of length $\leq m$
of elements 
of the form $m^r_N(\tilde v) i^r_N(\tilde w)$.
For each $m$ one has either\\  ${\rm dim} (\WWW^N_{m+1})>{\rm dim} (\WWW^N_m$) or 
$\WWW^N_m = \AAA^N_r$.\\  Let 
$D(N)=({\rm dim} \ S^N V)^2$. Since $D(N)\geq {\rm dim} \AAA^N_r$, one has $\WWW^N_{D(N)} =
\AAA^N_r$. 
In particular, $\WWW^{rd}_{D(rd)}$ contains elements of non-vanishing trace. Thus it 
suffices to take for 
$\nu(r,d)$ the least common multiple of all integers less or equal to $D(rd)$.

 $\hfill{\square}$

\section{ Proof of Lemma \ref{gbs}}

For each $x\in X$ and each non-zero $u\in (E_x\times F_x)^*$ we have to find a section
in $H^0(X,(E\otimes F)^{\otimes nr})$ which is not annihilated by $u^{\otimes nr}$.
Let 
$$I_E: H^0(X,S^rE)^{\otimes n}\rightarrow  H^0(X,E^{\otimes rn}),$$
$$I_F: H^0(X,S^rF)^{\otimes n}\rightarrow  H^0(X,F^{\otimes rn}),$$
$$J: H^0(X,E^{\otimes rn})\otimes H^0(X,F^{\otimes rn})\rightarrow H^0(X,(E\otimes F)^{\otimes rn})$$
be the canonical morphisms.
The permutation group $S(rn)$ acts on $E^{\otimes rn}$ in the standard way, which yields a map
$$\Sigma:S(rn)\times H^0(X,E^{\otimes rn}) \rightarrow  H^0(X,E^{\otimes rn}).$$
Let
$$\Phi:S(rn)\times H^0(X,S^rE)^{\otimes n}\otimes H^0(X,S^rF)^{\otimes n}\rightarrow H^0(X,(E\otimes F)^{\otimes nr})$$
be defined by $\Phi=J\circ(\Sigma\otimes Id)\circ(Id\times I_E\otimes I_F)$.
Though $\Phi$ factors through the coset map $S(rn)\rightarrow S(rn)/S(r)^{\times n}$, 
it yields sufficiently
many sections for our purpose.
For $x\in X$, $\sigma\in S(rn)$, $\tilde s_i\in H^0(X,S^rE)$, $\tilde t_i\in H^0(X,S^rF)$, 
there is an element
$$\gamma\in \Gamma_{E_x,F_x}(n,r,\tilde s_1(x),\ldots,\tilde s_n(x),\tilde t_1(x),
\ldots,\tilde t_n(x),u),$$ 
such that

\begin{equation} \label{phi}
\langle\Phi(\sigma,\tilde s_1\otimes \cdots \otimes\tilde s_n\otimes \tilde t_1\otimes 
\cdots \otimes\tilde t_n)_x,
 u^{\otimes nr}\rangle=\vert\gamma\vert.
\end{equation}

Conversely, for each $\gamma\in \Gamma(n,r,\tilde s_1(x),\ldots,\tilde s_N(x),\tilde t_1(x),
\ldots,\tilde t_n(x),u)$
one can find a permutation $\sigma\in S(rn)$ satisfying  the equation\eqref{phi}.

By Proposition \ref{nu} there exist sections $\tilde s_i$, $\tilde t_i$ and 
$$\gamma\in \Gamma(n,r,\tilde s_1(x),\ldots,\tilde s_N(x),\tilde t_1(x),\ldots,
\tilde t_n(x),u)$$ 
such that $\vert\gamma\vert\neq 0$. The image of the corresponding section 
$\Phi(\sigma,\tilde s_1\otimes \cdots \otimes\tilde s_n\otimes \tilde t_1\otimes \cdots 
\otimes\tilde t_n)$
in $H^0(X,S^{\otimes nr}(E\otimes F))$ yields a section in  $H^0(\PP (E\otimes F),
\OOO_{\PP (E\otimes F)}(nr))$
which generates the fibre over $(x,u^\perp)$.
 $\hfill{\square}$
 
 \begin{rem}
  The arguments in sections 1 and 2 remain valid when $\CC$ is replaced by any algebraically closed field of characteristic $0$ .
  Thus tensor products of semiample vector bundles on algebraic varieties over such a field are semiample, too.
 \end{rem}

\section {Direct sums and tensor products of $k$-ample vector bundles}
The following definition was introduced by Sommese \cite{S1}.

\begin{Def}  A line bundle $L$ on a compact complex space $X$ is $k$-ample if
 
1) L is semiample

2) the fibers of the corresponding morphism
$$\phi:X \to \PP H^0(X,L^r)^* $$
have dimensions less or equal to $k. $

A vector bundle $E$ is said to be $k$-ample if $\OOO_{\PP E}(1)$
is $k$-ample. 
\end{Def}

Note that $0$-ample is the same as ample, and $k$-ample implies $(k+1)$-ample.
Moreover $(\dim \PP E)$-ample is the same as semiample.

Sommese proved in (\cite{S1} p. 235, Corollary (1.10)) 
that direct sums and tensor products of $k$-ample vector bundles 
are $k$-ample, if these vector bundles are generated by sections. 
The aim of the present section is to remove the latter restriction, 
moreover for the product it is sufficient that one factor is k-ample 
and the other is semiample.

\begin{thm}\label{sum}

$E\oplus F$ is  $k$-ample if and only if $E,F$ are  $k$-ample.
\end{thm}

{\it Proof}: Since the quotient bundle of a $k$-ample bundle is $k$-ample, the 'only if' part 
is obvious. For the 'if' part we use criterion (1.7.3) of proposition (1.7) of Sommese 
\cite{S1}.
First note that $\OOO_{\PP E }(r_1)$ gbs (see Notation \ref{Notgbs}) and $\OOO_{\PP F}(r_2)$ gbs imply
that $\OOO_{\PP E }(r)$ and $\OOO_{\PP F }(r)$ are gbs whenever $r$ is a common multiple of $r_1,r_2$. 
By Lemma  \ref{sumgbs} this implies that $\OOO_{\PP (E\oplus F)}(r)$ is gbs. 
Now assume that there is a holomorphic finite to one map $\phi:Z\rightarrow X$
of a compact analytic space $Z$ to $X$, such that ${\rm dim} Z = k+1$ and that there is a surjective map 
$\phi^*(E\oplus F)\twoheadrightarrow Q$
with a trivial bundle $Q$. In particular there is a non-trivial section of
$\phi^*(E^*\oplus F^*)$. We may assume that the component $s:Z \rightarrow \phi^* E^*$
of this section is non-trivial.
Since $\OOO_{\PP E}(r)$ is gbs, for any $z\in Z$ such that $s(z)\neq 0$ there is a section
$\sigma: X\rightarrow S^rE$ such that $\langle \phi^*\sigma(z), s(z)^{\otimes r}\rangle\neq 0$.
In particular, $\langle \phi^*\sigma, s^r\rangle$ yields a non-trivial section
of the trivial line bundle over $Z$. Such a section must be constant, which implies
that $s$ cannot vanish anywhere. Consequently it yields a trivial quotient bundle of $\phi^*E$,
contrary to the assumption that $E$ is $k$-ample. $\hfill{\square}$

For the main theorem of this section we need two lemmata. We first reformulate  
criterion (1.7.4), proposition (1.7)  of Sommese \cite{S1}.

\begin{lem}\label{NG}
Let $E$ be semiample.
Then $E$ being $k$-ample is equivalent to the following condition:
Given any coherent sheaf $\GGG$ on $X$, there exists an integer $N(\GGG)$ 
such that for $n> N(\GGG)$ and $j>k$ one has
$$H^j(X,\GGG\otimes E^{\otimes n})=0.$$
\end{lem}

{\it Proof}: In its original form, criterion (1.7.4) uses the symmetric product $S^nE$
instead of
$E^{\otimes n}$. Since the symmetric product is a direct summand of $E^{\otimes n}$,
the new criterion implies the old. To show the converse, note that $E^{\otimes n}$
is isomorphic to a direct sum of direct summands of $S^n(E^{\oplus rk E})$. Moreover by theorem \ref{sum} 
$ E^{\oplus rk E}$ being $k$-ample is equivalent to $E$ being $k$-ample.

$\hfill{\square}$

\begin{lem}\label{filtration}
Let $F$ be a semiample vector bundle on a smooth complex space $X$. Then there are positive integers 
$r,{\tt L}$ and $n_l$ for $l=0,1,\ldots$, and graded locally free algebraic sheaves $\hat F_1,\ldots,\hat F_{\tt L}$ 
of finite rank with the following properties. The graded sheaf $SF=\oplus_{n\in\NN} S^nF$ has a resolution 
$$0 \rightarrow F_{\tt L} \rightarrow \ldots \rightarrow F_0 \rightarrow SF \rightarrow 0,$$
where each $F_i$ has a filtration
$$F_i=F_i^0\supset F_i^1\supset F_i^2\supset\ldots$$
compatible with the grading such that for $l=0,1,\ldots$
$$F_i^l/F_i^{l+1}\simeq \hat F_i[-lr]^{\oplus n_l},$$
where $[s]$ denotes the shift in degree by $s$.
\end{lem}

{\it Proof}: 
Choose $r$ such that  $\OOO_{\PP F}(r)$ is gbs. Let $\Gamma$ be the (isomorphic) image of 
$H^0(\PP F,\OOO_{\PP F}(r))$ in $H^0(X,S^rF)$. The elements of $S\Gamma$ act naturally on $SF$ by multiplication, in particular on $1 \in S^0F. $ 
Let $ \mathfrak{F}$ be the smallest subsheaf of $SF$ which contains the subsheaf $H^0(X,S^rF)1$ of $S^rF$ and admits 
the multiplicative action of $SF$. By Macaulay's theorem we have
$S^nF\subset \mathfrak{F}$ if $ n>N$, $N= (r-1) rk(F)$. 
By induction on $s$ one has $S^{sr+n}F\subset S^s\Gamma S^nF$ for $s\in\NN$ if $r+n>N$. 
The direct summands of $SF$ can be written in the form $S^{sr+n}F$ with $n\leq N$. Thus the multiplication map
$$S\Gamma \otimes (\oplus_{n=0}^N S^nF) \rightarrow SF$$ 
is surjective. This implies that $SF$ is a sheaf of finitely generated $(S\Gamma \otimes\OOO_X)$-modules.
Since $S\Gamma \otimes\OOO_X$ is a sheaf of regular local rings, $SF$ has local minimal resolutions
of finite length ${\tt L}$ which are unique up to isomorphism. Gluing the corresponding free $(S\Gamma \otimes\OOO_U)$-modules
over open sets $U\subset X$ by these isomorphisms, one obtains a resolution of $SF$ by sheaves of
locally free $(S\Gamma \otimes\OOO_X)$-modules $F_i$ of finite rank. Let 
$$F_i^l=  S^l\Gamma F_i$$
for $i=0,\ldots,\tt L$, and $\hat F_i = F_i/F_i^1$. This implies
$$F_i^l/F_i^{l+1}\simeq S^l\Gamma\otimes \hat F_i$$
for $l=0,1,\ldots$. Thus the statement of the proposition is true when $n_l$ is the dimension of the vector space $S^l\Gamma$.
$\hfill{\square}$

\begin{thm} \label{main}
If $E$ is a $k$-ample and $F$ a semiample vector bundle on a compact complex manifold
$X$, then $E\otimes F$  is $k$-ample.
\end{thm}

{\it Proof}: 
The sheaf  $(E\otimes F)^{\otimes n}$ is isomorphic to a direct sum of direct summands of
$E^{\otimes n} \otimes SF'$, where ${F'}=F^{\oplus rk F}$. Thus it suffices to show that 
$H^j(X,\GGG\otimes S{F'} \otimes E^{\otimes n})$ vanishes for $j>k$ and sufficiently
large $n$. The notations of lemmata \ref{NG} and \ref{filtration} will be adapted to $F'$ instead of $F$
by adding primes. Let $n> N(\GGG\otimes \hat {F'}_i)$ for $i=1,\ldots,{\tt L}'$. 
Then $H^j(X,\GGG\otimes {F'}_i^l/{F'}_i^{l+1}\otimes E^{\otimes n})=0$ for each $l$ and $j>k$.
Since the filtration of ${F'}_i$ is finite at each degree this implies
$H^j(X,\GGG\otimes {F'}_i\otimes E^{\otimes n})=0$ for $j>k$, $i=1,\ldots,{\tt L}'$,
and consequently $H^j(X,\GGG\otimes SF'\otimes E^{\otimes n})=0$.

$\hfill{\square}$

\begin{rem}
 Theorem \ref{main} implies that there is no need for the gbs conditions in 
 lemma (1.11.4) and proposition (1.13) of \cite{S1} when $X$ is smooth. The lemma and the  proposition can be strengthened as follows:
\end{rem}

\begin{lem}\label{Somm1}
 Let $E$ be a $k$-ample vector bundle on a compact complex manifold $X$, $Gr(s,E)$
 be the bundle of $s$-codimensional subspaces
of the fibres of $E$ and $\xi(s,E)$ the tautological line bundle of $Gr(s,E)$.
Then $\xi(s,E)$ is $k$-ample.
\end{lem}

In the proof of this lemma, $k$-ampleness of $E$ only enters through the $k$-ampleness of $\Lambda^sE$,
which is a quotient bundle of $E^{\otimes s}$.

\begin{prop}\label{Somm2} 
For $i=1,\ldots,n$ let $E_i$ be $k_i$-ample vector bundles of rank $r_i$ on 
a projective manifold $X$. 
$$H^p(X,\wedge^qT^*_X\otimes \wedge^{s_1}E_1\otimes \cdots \otimes \wedge^{s_n}E_n)=0$$
if $p+q>{\rm dim} X + \sum_i s_i(r_i-s_i) +\min \{ k_1, \ldots, k_n \}.$
\end{prop}

In the proof, $k$-ampleness of the $E_i$ only enters through the $k$-ampleness of 
$\wedge^{s_1}E_1\otimes \cdots \otimes \wedge^{s_n}E_n$, which is a quotient bundle of 
$E_1^{\otimes s_1}\otimes \cdots \otimes E_n^{\otimes s_n}$. 

See also \cite{LN}, where lemma 2.1 can easily be generalized to k-ample.

Finally we note two corollaries to theorems \ref{sum} and \ref{main}

\begin{corol}

If $E$ is a $k_1$-ample vector bundle and $F$ is $k_2$-ample vector bundle on a complex manifold, then 
$E\oplus F$ is $\max\{k_1,k_2\}$-ample and $E\otimes F$ is $\min \{ k_1,k_2 \}$-ample.
\end{corol}

\begin{corol}
One can extend the dominance ampleness theorem; theorem  3.7 in  \cite{dom} to an analogous
dominance $k$-ampleness theorem.
\end{corol}

Indeed, the proof of the dominance ampleness theorem uses no other property of ampleness but the fact that
direct summands, direct sums and tensors products of ample vector bundles are ample.

\end{document}